%% file: Perron_091205.tex
\input amstex
\let\myfrac=\frac
\input eplain
\let\frac=\myfrac
\input epsf

% Pour utiliser dvips, il faut employer :

% dvips -N0 -Z0 -K0 whatever.dvi -o whatever.ps

% Here we load the functions permitting us to use amsmath without using amsppt.

\loadeufm \loadmsam \loadmsbm
\message{symbol names}\UseAMSsymbols\message{,}

\font\myfontdefault=cmr10

% cmbx10 at 16pt%
\font\mytdmchapfont=cmb10 at 14pt
\font\mytdmheadfont=cmb10 at 10pt
\font\mytdmsubheadfont=cmr10

\magnification 1200
\newif\ifinappendices
\newif\ifundefinedreferences
\newif\ifchangedreferences
\newif\ifloadreferences
\newif\ifmakebiblio
\newif\ifmaketdm

\undefinedreferencestrue
\changedreferencesfalse

% With \loadreferencesfalse, the file makes its own counters and saves them to "references.tex". With \loadreferencestrue, the file loads the counters
% saved in "references.tex".

\loadreferencestrue
\makebibliofalse
\maketdmfalse

\def\headpenalty{-400}     % 400
\def\proclaimpenalty{-200} % 200

%%%%%%%%%%%%%%%%%%%%%%%%%%%%%%%%%%%%%%%%%%%%%%%%%%%%%%%%%%%%%%%%%%%%%%%%%%%%%%%%%%%%%%%%%%%%%%%%%%%%%%%%%%%%%%%%%%%%%%%
%
% Compteurs
%
%%%%%%%%%%%%%%%%%%%%%%%%%%%%%%%%%%%%%%%%%%%%%%%%%%%%%%%%%%%%%%%%%%%%%%%%%%%%%%%%%%%%%%%%%%%%%%%%%%%%%%%%%%%%%%%%%%%%%%%

\def\alphanum#1{\ifcase #1 _\or A\or B\or C\or D\or E\or F\or G\or H\or I\or J\or K\or L\or M\or N\or O\or P\or Q\or R\or S\or T\or U\or V\or W\or X\or Y\or Z\fi}
\def\gobbleeight#1#2#3#4#5#6#7#8{}

\newwrite\references
\newwrite\tdm
\newwrite\biblio

\newcount\chapno
\newcount\headno
\newcount\subheadno
\newcount\procno
\newcount\figno
\newcount\citationno

\def\setcatcodes{%
\catcode`\!=0 \catcode`\\=11}%

\ifloadreferences
    {\catcode`\@=11 \catcode`\_=11%
    \input references.tex %
    }%
\else
    \openout\references=references.tex
\fi

\newcount\newchapflag % Flag used for checking that we have a new chapter.
\newcount\showpagenumflag % Flag used for showing page numbers.

\global\chapno = -1 % We set it like this so that the introduction does not show up a number.
\global\citationno=0
\global\headno = 0
\global\subheadno = 0
\global\procno = 0
\global\figno = 0

\def\resetcounters{%
\global\headno = 0%
\global\subheadno = 0%
\global\procno = 0%
\global\figno = 0%
}

\global\newchapflag=0 % default - false.
\global\showpagenumflag=0 % default - false.

\def\chinfo{\ifinappendices\alphanum\chapno\else\the\chapno\fi}%
\def\headinfo{\ifinappendices\alphanum\headno\else\the\headno\fi}%
\def\subheadinfo{\headinfo.\the\subheadno}%      {\chinfo.\the\headno.\the\subheadno}%
\def\procinfo{\headinfo.\the\procno}%            {\chinfo.\the\headno.\the\procno}%
\def\figinfo{\the\figno}        %{\headinfo.\the\figno}%              {\chinfo.\the\headno.\the\figno}%
\def\citationinfo{\the\citationno}%
\def\nextheadno{\global\advance\headno by 1 \global\subheadno = 0 \global\procno = 0}% \headinfo\ - }
\def\nextsubheadno{\global\advance\subheadno by 1}% \subheadinfo\ - }
\def\nextprocno{\global\advance\procno by 1 \procinfo}
\def\nextfigno{\global\advance\figno by 1 \figinfo}

{\global\let\noe=\noexpand%
%
% Ici, je definit les macros qui me permettent de contruire des compteurs. Je vais changer provisoirement les catcodes afin de pouvoir bien ecrire le fichier
% contenant les references.
%
\catcode`\@=11%
\catcode`\_=11%
\setcatcodes%
!global!def!_@@internal@@makeref#1{%
!global!expandafter!def!csname #1ref!endcsname##1{%
!csname _@#1@##1!endcsname%
!expandafter!ifx!csname _@#1@##1!endcsname!relax%
    !write16{#1 ##1 not defined - run saving references}%
    !undefinedreferencestrue%
!fi}}%
!global!def!_@@internal@@makelabel#1{%
!global!expandafter!def!csname #1label!endcsname##1{%
!edef!temptoken{!csname #1info!endcsname}%
!ifloadreferences%
    !expandafter!ifx!csname _@#1@##1!endcsname!relax%
        !write16{#1 ##1 not hitherto defined - rerun saving references}%
        !changedreferencestrue%
    !else%
        !expandafter!ifx!csname _@#1@##1!endcsname!temptoken%
        !else
            !write16{#1 ##1 reference has changed - rerun saving references}%
            !changedreferencestrue%
        !fi%
    !fi%
!else%
    !expandafter!edef!csname _@#1@##1!endcsname{!temptoken}%
    !edef!textoutput{!write!references{\global\def\_@#1@##1{!temptoken}}}%
    !textoutput%
!fi}}%
!global!def!makecounter#1{!_@@internal@@makelabel{#1}!_@@internal@@makeref{#1}}%
!unsetcatcodes%
}
\makecounter{ch}%
\makecounter{head}%
\makecounter{subhead}%
\makecounter{proc}%
\makecounter{fig}%
\makecounter{citation}%
\def\newref#1#2{%
\def\temptext{#2}%
\edef\bibliotextoutput{\expandafter\gobbleeight\meaning\temptext}%
\global\advance\citationno by 1\citationlabel{#1}%
\ifmakebiblio%
    \edef\fileoutput{\write\biblio{\noindent\hbox to 0pt{\hss$[\the\citationno]$}\hskip 0.2em\bibliotextoutput\medskip}}%
    \fileoutput%
\fi}%
\def\cite#1{%
$[\citationref{#1}]$%
\ifmakebiblio%
    \edef\fileoutput{\write\biblio{#1}}%
    \fileoutput%
\fi%
}%
%
%%%%%%%%%%%%%%%%%%%%%%%%%%%%%%%%%%%%%%%%%%%%%%%%%%%%%%%%%%%%%%%%%%%%%%%%%%%%%%%%%%%%%%%%%%%%%%%%%%%%%%%%%%%%%%%%%%%%%%%
%
% Mise en page
%
%%%%%%%%%%%%%%%%%%%%%%%%%%%%%%%%%%%%%%%%%%%%%%%%%%%%%%%%%%%%%%%%%%%%%%%%%%%%%%%%%%%%%%%%%%%%%%%%%%%%%%%%%%%%%%%%%%%%%%%

\let\mypar=\par

% Here we define right justification.

\def\raggedleft{\leftskip=0pt plus 1fil \parfillskip=0pt}

% Lettrine - grace a Serroul

\font\lettrinefont=cmr10 at 28pt
\def\lettrine #1[#2][#3]#4%
{\hangafter -#1 \hangindent #2
\noindent\hskip -#2 \vtop to 0pt{
\kern #3 \hbox to #2 {\lettrinefont #4\hss}\vss}}

\font\mylettrinefont=cmr10 at 28pt
\def\mylettrine #1[#2][#3][#4]#5%
{\hangafter -#1 \hangindent #2
\noindent\hskip -#2 \vtop to 0pt{
\kern #3 \hbox to #2 {\mylettrinefont #5\hss}\vss}}

% Here we define the default page headers and footers.

\edef\Pagetitle={Blank}

\headline={\hfil\Pagetitle\hfil}
%{
%\ifnum\showpagenumflag=0
%    \hfil
%\else
%    \ifnum\newchapflag=0
%        \ifodd\pageno
%            \hfil\tenrm{\sl\Pagetitle}\ -\ \the\pageno
%        \else
%            \tenrm\hss\kern-1.2mm\the\pageno\ -\ {\sl\Pagetitle}\hfill
%        \fi
%    \else
%        \hfil
%    \fi
%\fi
%}

\footline={\hfil\myfontdefault\folio\hfil}
%{
%\ifnum\showpagenumflag=0
%    \hfil
%\else
%    \ifnum\newchapflag=0
%        \hfil
%    \else
%        \global\newchapflag=0
%        \ifodd\pageno
%            \hfil\tenrm\the\pageno
%        \else
%            \tenrm\the\pageno\hfil
%        \fi
%    \fi
%\fi
%}

\def\nextoddpage
{
\newpage%
\ifodd\pageno%
\else%
    \global\showpagenumflag = 0%
    \null%
    \vfil%
    \eject%
    \global\showpagenumflag = 1%
\fi%
}

% Ici on definit les titres des chapitres

\def\newchap#1#2%
{%
%
% First we reset the counters.
%
\global\advance\chapno by 1%
\resetcounters%
%
% Next we move on to the next page.
%
\newpage%
\ifodd\pageno%
\else%
    \global\showpagenumflag = 0%
    \null%
    \vfil%
    \eject%
    \global\showpagenumflag = 1%
\fi%
\global\newchapflag = 1%
\global\showpagenumflag = 1%
%
% Now we write the chapter heading.
%
{\font\chapfontA=cmsl10 at 30pt%
\font\chapfontB=cmsl10 at 25pt%
\null\vskip 5cm%
{\chapfontA\raggedleft\hfil%
{%
\ifnum\chapno=0
    \phantom{%
    \ifinappendices%
        Annexe \alphanum\chapno%
    \else%
        \the\chapno%
    \fi}%
\else%
    \ifinappendices%
        Annexe \alphanum\chapno%
    \else%
        \the\chapno%
    \fi%
\fi%
}%
\par}%
\vskip 2cm%
{\chapfontB\raggedleft%
\lineskiplimit=0pt%
\lineskip=0.8ex%
\hfil #1\par}%
\vskip 2cm%
}%
\edef\Pagetitle{#2}%
%
% Finally, we write the information into the tdm file.
%
\ifmaketdm%
    \def\temp{#2}%
    \def\tempbis{\nobreak}%
    \edef\chaptitle{\expandafter\gobbleeight\meaning\temp}%
    \edef\mynobreak{\expandafter\gobbleeight\meaning\tempbis}%
    \edef\textoutput{\write\tdm{\bigskip{\noexpand\mytdmchapfont\noindent\chinfo\ - \chaptitle\hfill\noexpand\folio}\par\mynobreak}}%
\fi%
\textoutput%
}

% Ici on definit les sous titres.

\def\newhead#1%
{%
\ifhmode%
    \mypar%
\fi%
\ifnum\headno=0%
\ifinappendices
    \nobreak\vskip -\lastskip%
    \nobreak\vskip .5cm%
\fi
\else%
    \nobreak\vskip -\lastskip%
    \nobreak\vskip .5cm%
\fi%
\nextheadno%
\ifmaketdm%
    \def\temp{#1}%
    \edef\sectiontitle{\expandafter\gobbleeight\meaning\temp}%
    \edef\textoutput{\write\tdm{\noindent{\noexpand\mytdmheadfont\quad\headinfo\ - \sectiontitle\hfill\noexpand\folio}\par}}%
    \textoutput%
\fi%
\font\headfontA=cmbx10 at 14pt%
{\headfontA\noindent\headinfo\ - #1.\hfil}%
\nobreak\vskip .5cm%
}%

% Ici on definit les sous-sous titres.

\def\newsubhead#1%
{%
\ifhmode%
    \mypar%
\fi%
\ifnum\subheadno=0%
\else%
    \penalty\headpenalty\vskip .4cm%
\fi%
\nextsubheadno%
\ifmaketdm%
    \def\temp{#1}%
    \edef\subsectiontitle{\expandafter\gobbleeight\meaning\temp}%
    \edef\textoutput{\write\tdm{\noindent{\noexpand\mytdmsubheadfont\quad\quad\subheadinfo\ - \subsectiontitle\hfill\noexpand\folio}\par}}%
    \textoutput%
\fi%
\font\subheadfontA=cmsl10 at 12pt% cmbxsl10 at 12pt% cmbxti110 at 12pt%
{\subheadfontA\noindent\subheadinfo\ #1.\hfil}%
\nobreak\vskip .25cm %
}%

%%%%%%%%%%%%%%%%%%%%%%%%%%%%%%%%%%%%%%%%%%%%%%%%%%%%%%%%%%%%%%%%%%%%%%%%%%%%%%%%%%%%%%%%%%%%%%%%%%%%%%%%%%%%%%%%%%%%%%%
%
% Commandes et Symboles
%
%%%%%%%%%%%%%%%%%%%%%%%%%%%%%%%%%%%%%%%%%%%%%%%%%%%%%%%%%%%%%%%%%%%%%%%%%%%%%%%%%%%%%%%%%%%%%%%%%%%%%%%%%%%%%%%%%%%%%%%

% Ici on definit la famille mathroman...

\font\mathromanten=cmr10
\font\mathromanseven=cmr7
\font\mathromanfive=cmr5
\newfam\mathromanfam
\textfont\mathromanfam=\mathromanten
\scriptfont\mathromanfam=\mathromanseven
\scriptscriptfont\mathromanfam=\mathromanfive
\def\mathroman{\fam\mathromanfam}

% Ici on definit la famille mathsf qui n'est pas definie dans Plain.

\font\sf=cmss12

\font\sansseriften=cmss10
\font\sansserifseven=cmss7
\font\sansseriffive=cmss5
\newfam\sansseriffam
\textfont\sansseriffam=\sansseriften
\scriptfont\sansseriffam=\sansserifseven
\scriptscriptfont\sansseriffam=\sansseriffive
\def\mathsf{\fam\sansseriffam}

% Ici on definit la famille mathbf qui n'est pas definie dans Plain.

\font\bftwelve=cmb12

\font\boldten=cmb10
\font\boldseven=cmb7
\font\boldfive=cmb5
\newfam\mathboldfam
\textfont\mathboldfam=\boldten
\scriptfont\mathboldfam=\boldseven
\scriptscriptfont\mathboldfam=\boldfive
\def\mathbf{\fam\mathboldfam}

% Ici on definit \mathi et \mathj. En fait, cela n'est pas necessaire, car ces symboles existent deja sous les noms respectivement de \imath et \jmath. Mais,
% ca reste une exercice interessant.

\font\mycmmiten=cmmi10
\font\mycmmiseven=cmmi7
\font\mycmmifive=cmmi5
\newfam\mycmmifam
\textfont\mycmmifam=\mycmmiten
\scriptfont\mycmmifam=\mycmmiseven
\scriptscriptfont\mycmmifam=\mycmmifive

\def\hexa#1{\ifcase #1 0\or 1\or 2\or 3\or 4\or 5\or 6\or 7\or 8\or 9\or A\or B\or C\or D\or E\or F\fi}
\mathchardef\mathi="7\hexa\mycmmifam7B
\mathchardef\mathj="7\hexa\mycmmifam7C

% Ici je definit \mybeth, \mygimmel et \mydaleth

\font\mymsbmten=msbm10 at 8pt%8.5pt
\font\mymsbmseven=msbm7 at 5.6pt%6pt
\font\mymsbmfive=msbm5 at 4pt%4.25pt
\newfam\mymsbmfam
\textfont\mymsbmfam=\mymsbmten
\scriptfont\mymsbmfam=\mymsbmseven
\scriptscriptfont\mymsbmfam=\mymsbmfive

\mathchardef\mybeth="7\hexa\mymsbmfam69
\mathchardef\mygimmel="7\hexa\mymsbmfam6A
\mathchardef\mydaleth="7\hexa\mymsbmfam6B

% Ici je definit la fonction qui insere des figures.

\def\placelabel[#1][#2]#3{{%
\setbox10=\hbox{\raise #2cm \hbox{\hskip #1cm #3}}%
\ht10=0pt%
\dp10=0pt%
\wd10=0pt%
\box10}}%

% Ici je definit les fonctions \proclaim et \endproclaim. Les premieres versions sont plus jolies mais je ne sais pas comment les couper verticalement, et il y
% a alors des problemes de mise en page. Les deuxiemes versions sont plus simples, mais il n'y a pas de problem de mise en page.

\newif\ifinproclaim%
\global\inproclaimfalse%
\def\proclaim#1{%
\medskip%
%
% J'insere ici un bgroup pour que les variables que je definis ici soit uniquement locales.
%
\bgroup%
\inproclaimtrue%
\setbox10=\vbox\bgroup\leftskip=0.8em\noindent{\bftwelve #1}\sf%
}

\def\endproclaim{%
\egroup%
\setbox11=\vtop{\noindent\vrule height \ht10 depth \dp10 width 0.1em}%
\wd11=0pt%
\setbox12=\hbox{\copy11\kern 0.3em\copy11\kern 0.3em}%
\wd12=0pt%
\setbox13=\hbox{\noindent\box12\box10}%
\noindent\unhbox13%
\egroup%
\medskip\ignorespaces%
}

\def\proclaim#1{%
\medskip%
\bgroup%
\inproclaimtrue%
\noindent{\bftwelve #1}%
\nobreak\medskip%
\sf%
}

\def\endproclaim{%
\mypar\egroup\penalty\proclaimpenalty\medskip\ignorespaces%
}

\def\noskipproclaim#1{%
\medskip%
\bgroup%
\inproclaimtrue%
\noindent{\bf #1}\nobreak\sl%
}

\def\endnoskipproclaim{%
\mypar\egroup\penalty\proclaimpenalty\medskip\ignorespaces%
}

% Ici on definit quelques autres commandes dont on aura besoin, en particulier, le colon francais.

\def\proof{{\noindent\bf Proof:\ }}

\def\remark{{\noindent\sl Remark:\ }}

\def\msf#1{{\mathsf #1}}

\def\qed{~$\square$}
\def\munion{\mathop{\cup}}
\def\minter{\mathop{\cap}}
\def\myitem#1{%
%\ifinproclaim%
%    \item{#1}%
%\else%
    \noindent\hbox to .5cm{\hfill#1\hss}%\kern 0.1em}%
}
%\fi}%

\catcode`\@=11
\def\Eqalign#1{\null\,\vcenter{\openup\jot\m@th\ialign{%
\strut\hfil$\displaystyle{##}$&$\displaystyle{{}##}$\hfil%
&&\quad\strut\hfil$\displaystyle{##}$&$\displaystyle{{}##}$%
\hfil\crcr #1\crcr}}\,}
\catcode`\@=12

\def\makeop#1{%
\global\expandafter\def\csname op#1\endcsname{{\mathroman #1}}}%

\def\makeopsmall#1{%
\global\expandafter\def\csname op#1\endcsname{{\mathroman{\lowercase{#1}}}}}%

\makeopsmall{ArcTan}%
\makeopsmall{ArcCos}%
\makeop{Arg}%
\makeop{Det}%
\makeop{Log}%
\makeop{Re}%
\makeop{Im}%
\makeop{Dim}%
\makeopsmall{Tan}%
\makeop{Ker}%
\makeopsmall{Cos}%
\makeopsmall{Sin}%
\makeop{Exp}%
\makeopsmall{Tanh}%
\makeop{Tr}%
\makeop{End}%
\makeop{Long}%
\makeop{Ch}%
\makeop{Exp}%
\makeop{Eval}%
\makeop{Lift}%
\makeop{Int}%
\makeop{Ext}%
\makeop{Aire}%
\makeop{Im}%
\makeop{Conf}%
\makeop{Exp}%
\makeop{Mod}%
\makeop{Log}%
\makeop{Sgn}%
\makeop{Ext}%
\makeop{Int}%
\makeop{Ln}%
\makeop{Dist}%
\makeop{Aut}%
\makeop{Id}%
\makeop{GL}%
\makeop{SO}%
\makeop{Homeo}%
\makeop{Vol}%
\makeop{Ric}%
\makeop{Hess}%
\makeop{Euc}%
\makeop{Isom}%
\makeop{Max}%
\makeop{SW}%
\makeop{SL}%
\makeop{Long}%
\makeop{Fix}%
\makeop{Wind}%
\makeop{Diag}%
\makeop{dVol}%
\makeop{Symm}%
\makeop{Ad}%
\makeop{Diam}%
\makeop{loc}%
\makeopsmall{Sinh}%
\makeopsmall{Cosh}%
\makeop{Len}%
\makeop{Length}%
\makeop{Conv}%
\makeop{Min}%
\makeop{Area}%
\font\mycirclefont=cmsy7
\def\textcircle{{\raise 0.3ex \hbox{\mycirclefont\char'015}}}

\let\emph=\bf

\hyphenation{quasi-con-formal}

%%%%%%%%%%%%%%%%%%%%%%%%%%%%%%%%%%%%%%%%%%%%%%%%%%%%%%%%%%%%%%%%%%%%%%%%%%%%%%%%%%%%%%%%%%%%%%%%%%%%%%%%%%%%%%%%%%%%%%%
%
% Citations
%
%%%%%%%%%%%%%%%%%%%%%%%%%%%%%%%%%%%%%%%%%%%%%%%%%%%%%%%%%%%%%%%%%%%%%%%%%%%%%%%%%%%%%%%%%%%%%%%%%%%%%%%%%%%%%%%%%%%%%%%

\ifmakebiblio%
    \openout\biblio=biblio.tex %
    {%
        \edef\fileoutput{\write\biblio{\bgroup\leftskip=2em}}%
        \fileoutput
    }%
\fi%

\newref{HarveyLawsonI}{Harvey F. R., Lawson H. B. Jr., Dirichlet Duality and the Nonlinear Dirichlet Problem, {\sl Comm. Pure Appl. Math.} {\bf 62} (2009), no. 3, 396--443}
\newref{HarveyLawsonII}{}
\newref{CaffNirSprV}{Caffarelli L., Nirenberg L., Spruck J., Nonlinear second-order elliptic equations. V. The Dirichlet problem for Weingarten hypersurfaces, {\sl Comm. Pure Appl. Math.} {\bf 41} (1988), no. 1, 47--70}
\newref{CrandIshiiLions}{Crandall M. G., Ishii H., Lions P.-L., User's guide to viscosity solutions of second order partial differential equations, {\sl Bull. Amer. Math. Soc.} {\bf 27} (1992), no. 1, 1--67}
\newref{GuanSpruckI}{Guan B., Spruck J., The existence of hypersurfaces of constant Gauss curvature with prescribed boundary, {\sl J. Differential Geom.} {\bf 62} (2002), no. 2, 259--287}
\newref{GuanSpruckII}{}
\newref{SmiCGC}{Smith G., Constant Gaussian Curvature Hypersurfaces in Hadamard Manifold}
\newref{SmiSLC}{Smith G., Special Lagrangian curvature, arXiv:math/0506230}
\newref{SmiNLD}{Smith G., The Non-Linear Dirichlet Problem in Hadamard Manifolds, arXiv:0908.3590}

\ifmakebiblio%
    {\edef\fileoutput{\write\biblio{\egroup}}%
    \fileoutput}%
\fi%

%%%%%%%%%%%%%%%%%%%%%%%%%%%%%%%%%%%%%%%%%%%%%%%%%%%%%%%%%%%%%%%%%%%%%%%%%%%%%%%%%%%%%%%%%%%%%%%%%%%%%%%%%%%%%%%%%%%%%%%
%
% The Paper
%
%%%%%%%%%%%%%%%%%%%%%%%%%%%%%%%%%%%%%%%%%%%%%%%%%%%%%%%%%%%%%%%%%%%%%%%%%%%%%%%%%%%%%%%%%%%%%%%%%%%%%%%%%%%%%%%%%%%%%%%
%
\document
\myfontdefault
\global\chapno=1
\global\showpagenumflag=1
\def\Pagetitle{}
\null
\vfill
\def\centre{\rightskip=0pt plus 1fil \leftskip=0pt plus 1fil \spaceskip=.3333em \xspaceskip=.5em \parfillskip=0em \parindent=0em}%
\def\textmonth#1{\ifcase#1\or January\or Febuary\or March\or April\or May\or June\or July\or August\or September\or October\or November\or December\fi}
\font\abstracttitlefont=cmr10 at 14pt
{\abstracttitlefont\centre The Perron Method and the Non-Linear Plateau Problem\par}
\bigskip
{\centre \the\day\ \textmonth\month\ \the\year\par}
\bigskip
{\centre{\bf Andrew Clarke}\par}
\bigskip
{\centre Laboratoire de Math\'ematiques Jean Leray,\par
2, rue de la Houssini\`ere - BP 92208,\par
44322 Nantes Cedex 3,\par
FRANCE\par}
\bigskip
{\centre{\bf Graham Smith}\par}
\bigskip
{\centre Departament de Matem\`atiques,\par
Facultat de Ci\`encies, Edifici C,\par
Universitat Aut\`onoma de Barcelona,\par
08193 Bellaterra,\par
Barcelona,\par
SPAIN\par}
\bigskip
\noindent{\emph Abstract:\ }We describe a novel technique for solving the Plateau problem for constant curvature hypersurfaces based on recent work of Harvey and Lawson. This is illustrated by an existence theorem for hypersurfaces of constant Gaussian curvature in $\Bbb{R}^{n+1}$.
%
% 58E12 : Global analysis, analysis on manifolds; applications to minimal surfaces.
% 35J25 : Boundary value problems for second order elliptic equations
% 35J60 : Nonlinear elliptic equations
% 53A10 : Minimal surfaces, surfaces with prescribed mean curvature
% 53C21 : Methods of Riemannian geometry, including PDE methods; curvature restrictions
% 53C42 : Immersions (minimal, prescribed curvature, tight, etc.)
%
\par 
\vfill
\nextoddpage
\global\pageno=1
\def\Pagetitle{\sl The Perron Method and the Non-Linear Plateau Problem}
\newhead{Introduction}
\noindent In this paper we describe a novel technique for constructing solutions to the Plateau problem for convex hypersurfaces of constant Gaussian curvature, which we illustrate through the proof of the following theorem:
\proclaim{Theorem \nextprocno}
\noindent Let $\hat{K}\subseteq\Bbb{R}^{n+1}$ be a compact, strictly convex set with smooth boundary. Let $\Omega$ be an open subset of $\partial\hat{K}$ with (non-trivial) smooth boundary. Suppose there exists $k>0$ such that the Gaussian curvature of $\hat{K}$ is everywhere at least $k$. Then, for all $t\in]0,k]$, there exists a convex subset $K_t\subseteq\hat{K}$ such that:
\medskip
\myitem{(i)} $K_t\minter\partial\hat{K}=\Omega^c$; and
\medskip
\myitem{(ii)} the boundary of $K_t$ is smooth in the interior of $K$ and is of constant Gaussian curvature equal to $t$.
\endproclaim
\proclabel{TheoremExistence}
\remark This follows directly from Lemma \procref{LemmaExistence} and Theorem $5.1$ of \cite{GuanSpruckI}.
\medskip
\noindent Hypersurfaces of constant Gaussian curvature are interesting objects of study for various reasons. When $n=2$, the Gaussian curvature is (more or less) equivalent to the intrinsic curvature of the surface, which only depends on one variable. In higher dimensions, although no such relation exists, Gaussian curvature continues to provide relatively simple PDEs that make it a good model for the study of more general non-linear notions of curvature. A tremendous literature exists devoted to the study of this problem, of which the most significant results are perhaps \cite{CaffNirSprV} of Caffarelli, Nirenberg and Spruck and \cite{GuanSpruckI} of Guan and Spruck. We refer the reader to the introduction of the paper \cite{SmiCGC} by the second author for a broader overview.
\medskip
\noindent The novel technique that we introduce is a version of the Perron Method recently developed by Harvey and Lawson in \cite{HarveyLawsonI} and \cite{HarveyLawsonII}. This yields convex sets whose boundaries are of constant Gaussian curvature in the viscosity sense (c.f. \cite{CrandIshiiLions}). We then show that these hypersurfaces are smooth away from their boundaries by appealling to Theorem $5.1$ of \cite{GuanSpruckI}. The elementary nature of this proof as well as the remarkable generality of the results of \cite{HarveyLawsonI} and \cite{HarveyLawsonII} hint at potential generalisations. Indeed, it can easily be extended to yield (not necessarily unique) viscosity solutions in any Hadamard manifold, and the regularity result of \cite{GuanSpruckI} can then be applied whenever the ambient manifold is also affine flat (c.f. \cite{SmiCGC}). This is, in particular, the case for hyperbolic space. 
\medskip
\noindent This approach can also be extended to treat other notions of curvature. In fact, it is currently applicable to any notion of curvature which constitutes a ``convex condition'' in the sense of Section \headref{HeadPSets}. In particular, this includes special Lagrangian curvature, which has been studied extensively by the second author in \cite{SmiSLC} and \cite{SmiNLD}.
\medskip
\noindent An important aspect of this technique that departs from the approach of Harvey and Lawson is its dependance on convex sets. We have introduced this essentially in order to make the problem more tractable, and it does so in two ways. The first is by eliminating complicated geometric considerations such as self intersections. This is not an issue in \cite{HarveyLawsonI} and \cite{HarveyLawsonII}, since Harvey and Lawson are working there with functions, where the geometry is constant, as it were. The second is by allowing us to use the $C^{0,1}$ regularity properties of convex sets, and thus neatly sidestep the problem of proving regularity, which is the hardest step in Harvey and Lawson's proof and in the study of viscosity solutions in general. Naturally, however, this dependance on convex sets is very restrictive, and excludes large families of interesting curvature functions (see, for example, \cite{CaffNirSprV} or \cite{GuanSpruckII}). We nonetheless expect appropriate modifications to yield stronger results in the near future.
\medskip
\noindent This paper is structured as follows:
\medskip
\myitem{(i)} in Section $2$, we define the basic notions used throughout this paper, recalling the definition of Dirichlet set as introduced by Harvey and Lawson in \cite{HarveyLawsonI};
\medskip
\myitem{(ii)} in Section $3$, we define what it means for a convex set to be of type $F$. We show that this constitutes a ``Perron system'' in the sense that it satisfies the basic axioms required for the Perron Method; and
\medskip
\myitem{(iii)} in Section $4$, we apply the Perron Method to obtain viscosity solutions to the Plateau problem, and, appealing to Theorem $5.1$ of \cite{GuanSpruckI}, this proves Theorem \procref{TheoremExistence}.
\medskip
\noindent The second author would like to thank Prof. Lawson for bringing \cite{HarveyLawsonII} to his attention. 
\goodbreak
\newhead{Dirichlet Sets}
\noindent Let $\opSymm(\Bbb{R}^n)$ denote the space of symmetric matrices over $\Bbb{R}^n$. We define $P\subseteq\opSymm(\Bbb{R}^n)$ to be the set of all symmetric, non-negative semi-definite matrices. Thus $A$ is an element of $P$ if and only for all $x\in\Bbb{R}^n$:
\headlabel{HeadPSets}
$$
\langle Ax,x\rangle\geqslant 0.
$$
\noindent Trivially, $P$ is a closed convex cone. Let $F$ be a closed subset of $\opSymm(\Bbb{R}^n)$. Following \cite{HarveyLawsonI}, we will say that $F$ is a Dirichlet set if and only if:
$$
F+P\subseteq F.
$$
\noindent Moreover, we will say that $F$ is invariant if and only it is preserved by conjugation by matrices in $O(n)$. In other words, $F$ is invariant if and only if, for all $A\in F$, and for all $M\in O(n)$:
$$
M^tAM \in F.
$$
\noindent Finally, we will say that $F$ defines a convex condition if and only if:
$$
F\subseteq P.
$$
\noindent In this paper, we are interested in invariant Dirichlet sets which define convex conditions.
\medskip
{\sl\noindent Example:\ } For $k>0$ define $F_k$ by:
$$
F_k = \left\{A\in P\text{ s.t. }\opDet(A)\geqslant k\right\}.
$$
\noindent It is easily verified that $F_k$ is an invariant Dirichlet set which trivially also defines a convex condition. In fact, $F_k$ is also convex with smooth boundary.\qed
\goodbreak
\newhead{The Perron System}
\noindent Let $\msf{N}$ be a unit normal vector field over a smooth hypersurface, $\Sigma$. In the sequel, we adopt the convention whereby the shape operator, $A$, of $\Sigma$ satisfies:
$$
A\cdot X = \nabla_X\msf{N}.
$$
\noindent Let $F$ be an invariant Dirichlet set. Let $X\subseteq\Bbb{R}^{n+1}$ be a compact set. We say that $X$ is of type $F$ if and only if, for all $p\in\partial X$, if $\Omega$ is an open subset of $X^o$ (the interior of $X$) such that:
\medskip
\myitem{(i)} $\partial\Omega$ is smooth; and
\medskip
\myitem{(ii)} $p\in\partial\Omega$,
\medskip
\noindent then the shape operator of $\partial\Omega$ at $p$ with respect to the outward pointing normal is conjugate to a matrix in $F$.
\medskip
\remark The shape operator of $\partial\Omega$ is conjugate to a matrix in $F$ if and only if its matrix with respect to an orthonormal basis for $T\partial\Omega$ lies in $F$. Since $F$ is $O(n)$-invariant, this does not depend on the orthonormal basis chosen.
\medskip
\remark Observe that if $F$ is not a Dirichlet set, then this definition is essentially empty.
\medskip
\noindent Following \cite{HarveyLawsonI} and \cite{HarveyLawsonII}, we obtain the following characterisation of subsets of type $F$:
\proclaim{Lemma \nextprocno}
\noindent Let $X\subset\Bbb{R}^{n+1}$ be a compact set. $X$ is not of type $F$ if and only if there exists $p\in\partial X$, $r>0$ and $f:B_r(p)\rightarrow\Bbb{R}$ such that:
\medskip
\myitem{(i)} $f(p)=0$;
\medskip
\myitem{(ii)} $f^{-1}(]-\infty,0])\subseteq X\minter B_r(p)$;
\medskip
\myitem{(iii)} $(\nabla f)(q)\neq 0$ for all $q\in B_r(p)$; and
\medskip
\myitem{(iv)} $\opHess(f)|_{(\nabla f)^\perp}(q)$ is conjugate to an element of $\|\nabla f(q)\|F^c$ for all $q\in B_r(p)$.
\endproclaim
\proclabel{LemmaCharacterisation}
\proof We recall that if $f$ is smooth and if $\nabla f\neq 0$ at a point $q$, then $\nabla f$ is colinear with the normal vector to the level set of $f$ passing through $q$, which we denote by $\Sigma_q$. Moreover, if $A_q$ is the shape operator of $\Sigma_q$ with respect to the normal pointing in the same direction as $\nabla f$, then:
$$
A_q = \|(\nabla f)(q)\|^{-1}\opHess(f)|_{(\nabla f)^\perp}(q).
$$ 
\noindent The result follows directly from these relations and the contrapositive of the definition of being of type $F$.\qed
\medskip
\noindent In the case of smooth boundary, we obtain:
\proclaim{Lemma \nextprocno}
\noindent Let $X\subseteq\Bbb{R}^{n+1}$ be a compact set. Suppose that $\partial X$ is smooth, then $X$ is of type $F$ if and only if the shape operator of $\partial X$ with respect to the outward pointing normal is conjugate to an element of $F$ at every point of $\partial X$.
\endproclaim
\noindent Now suppose, moreover, that $F$ defines a convex condition. The Perron method is based on the following result:
\proclaim{Lemma \nextprocno}
\noindent Let $\Cal{F}$ be a family of compact, convex sets of type $F$. Let $X$ be the intersection of all members of $\Cal{F}$. Then $X$ is also a compact, convex set of type $F$.
\endproclaim
\proclabel{LemmaIntersectionIsOfTypeF}
\proof $X$ is trivially compact and convex. Suppose that $X$ is not of type $F$. By Lemma \procref{LemmaCharacterisation}, there exists $p\in\partial X$, $r>0$ and $f:B_r(p)\rightarrow\Bbb{R}$ such that:
\medskip
\myitem{(i)} $f(p)=0$;
\medskip
\myitem{(ii)} $f^{-1}(]-\infty,0])\subseteq X\minter B_r(p)$; 
\medskip
\myitem{(iii)} $(\nabla f)(q)\neq 0$ for all $q\in B_r(p)$; and
\medskip
\myitem{(iv)} $\opHess(f)|_{(\nabla f)^\perp}(q)$ is conjugate to an element of $\|\nabla f(q)\|F^c$ for all $q\in B_r(p)$.
\medskip
\noindent Moreover, since $F^c$ is open, by reducing $r$ if necessary, $f$ may be chosen such that there exists $\epsilon>0$ such that, for all $q\in X^c\minter \partial B_r(p)$:
$$
f(q)\geqslant \epsilon.
$$
\noindent Choose $q\in B_r(p)\minter X^c$ such that $f(q)\leqslant\epsilon/2$. There exists $Y\in\Cal{F}$ such that $q\notin Y$. However, since $X\subseteq Y$:
$$
f^{-1}(]-\infty,0])\subseteq Y.
$$
\noindent Let $p'\in B_r(p)$ be the point in the closure of $Y^c\minter B_r(p)$ realising the infimum of $f$ over this set, and let $\delta$ be the value of this infimum. Trivially $0\leqslant\delta\leqslant\epsilon/2$. Since, for all $q\in Y\minter\partial B_r(p)$, $f(q)\geqslant \epsilon$, $p'$ is an interior point of $B_r(p)$, so there exists $r'>0$ such that:
$$
B_{r'}(p')\subseteq B_r(p).
$$
\noindent Defining $f':B_{r'}(p)\rightarrow\Bbb{R}$ by $f'=f-\delta$, we deduce by Lemma \procref{LemmaCharacterisation} that $Y$ is not of type $F$. This contradicts the hypothesis on $\Cal{F}$, and the result follows.\qed
\goodbreak
\newhead{Duality and the Viscosity Solution}
\noindent Let $\hat{K}\subseteq\Bbb{R}^{n+1}$ be a compact, strictly convex subset with smooth boundary of Gaussian curvature at least $k>0$. Let $\Omega\subseteq\partial K$ be an open subset with smooth boundary. Let $K_0$ be the convex hull of $\Omega^c$. Observe that, since $\Omega$ has smooth boundary, $K_0$ has non-trivial interior.
\medskip
\noindent Let $X\subseteq\Bbb{R}^{n+1}$ be a compact set. We say that $X$ is of type $F'$ if and only if, for all $p\in\partial X$, if $\Omega$ is an open subset of $X^c$ such that:
\medskip
\myitem{(i)} $\partial\Omega$ is smooth; and
\medskip
\myitem{(ii)} $p\in\partial\Omega$,
\medskip
\noindent then the shape operator of $\partial\Omega$ at $p$ with respect to the inward pointing normal is conjugate to a matrix in $\overline{F^c}$.
\medskip
\remark In the language of \cite{HarveyLawsonI} and \cite{HarveyLawsonII}, this is the dual property to the property of being of type $F$. The duality becomes evident when we observe that $\tilde{F}:=-\overline{F^c}$ is also an invariant Dirichlet set (c.f. \cite{HarveyLawsonI}), and, when $X$ is the closure of its interior, $X$ is of type $F'$ if and only if $\overline{X^c}$ is of type $\tilde{F}$.
\medskip
\noindent For all $t>0$, we define $F_t\subseteq\opSymm(\Bbb{R}^n)$ by:
$$
F_t = \left\{A\in P\text{ s.t. }\opDet(A)\geqslant t\right\}.
$$
\noindent As discussed in section \headref{HeadPSets}, $F_t$ is an invariant Dirichlet set which defines a convex condition.
\proclaim{Lemma \nextprocno}
\noindent For all $t\in]0,k]$, there exists a compact, convex subset $K_t$ of $\hat{K}$ such that:
\medskip
\myitem{(i)} $K_0\subseteq K_t$; 
\medskip
\myitem{(ii)} $K_t\minter\partial\hat{K}=\Omega^c$; and
\medskip
\myitem{(iii)} $K_t$ is of type $F_t$ and of type $F_t'$ over the interior of $\hat{K}$.
\endproclaim
\proclabel{LemmaExistence}
\remark Thus, for all $t$, $K_t$ has constant Gaussian curvature in the viscosity sense (c.f. \cite{CrandIshiiLions}).
\medskip
\proof Choose $t\in]0,k]$. Let $\Cal{F}$ denote the set of all convex subsets of $\hat{K}$ which contain $K_0$ and which are of type $F_t$. $\Cal{F}$ is non-empty since $\hat{K}\in\Cal{F}$. Let $K_t$ be the intersection of all members of $\Cal{F}$. Trivially, $K_0\subseteq K_t$. Since $\hat{K}$ is strictly of type $F_t$, $K_t\minter\partial\hat{K}=\Omega^c$. Finally, by Lemma \procref{LemmaIntersectionIsOfTypeF}, $K_t$ is of type $F_t$ over the interior of $\hat{K}$. It thus remains to prove that $K_t$ is of type $F_t'$ over the interior of $\hat{K}$.
\medskip
\noindent Suppose the contrary. Observe that, since $K_0$ has non-trivial interior, so does $K_t$. Moreover, $K_t\in\Cal{F}$. By Lemma \procref{LemmaCharacterisation} (reversing orientation), there exists $p\in\partial K_t\minter \hat{K}$, $r>0$ and a smooth function $f:B_r(p)\rightarrow\Bbb{R}$ such that:
\medskip
\myitem{(i)} $f(p)=0$;
\medskip
\myitem{(ii)} $K_t\minter B_r(p)\subseteq f^{-1}(]-\infty,0])$;
\medskip
\myitem{(iii)} $(\nabla f)(q)\neq 0$ for all $q\in B_r(p)$; and
\medskip
\myitem{(iv)} $\opHess(f)|_{(\nabla f)^\perp}(q)$ is conjugate to an element of $\|(\nabla f)(q)\|F^o_t$ for all $q\in B_r(p)$, where $F^o_t$ is the interior of $F_t$.
\medskip
\noindent Moreover, since $F^o_t$ is open, by reducing $r$ if necessary, $f$ may be chosen such that there exists $\epsilon>0$ such that, for all $q\in K_t\minter \partial B_r(p)$:
$$
f(q)\leqslant -\epsilon.
$$
\noindent For all $t\in]-\epsilon,0]$, define $\Sigma_t$ by:
$$
\Sigma_t = f^{-1}(\left\{t\right\}).
$$
\noindent For all such $t$, $\partial\Sigma_t$, which lies in $\partial B_r(p)$ is a subset of $K_t^c$ and therefore also of $K_0^c$. Moreover, for all such $t$, $\Sigma_t$ is strictly convex over its interior. 
We recall that, since $K_0$ is a convex hull of a subset of $\partial\hat{K}$, $\partial K_0$ is locally ruled throughout the interior of $\hat{K}$. In other words, for all $p\in\partial K_0$ lying in the interior of $\hat{K}$, there exists a straight line segment, $\Gamma$, containing $p$ in its interior and which also lies in $\partial K_0$. Thus, since $\Sigma_0$ lies in the closure of the complement of $K_0$, we deduce by the geometric maximum principal that so does $\Sigma_t$ for all $t\in]-\epsilon,0]$. In particular, if we define $K_t'$ by:
$$
K_t' = (K_t\minter f^{-1}(]-\infty,-\epsilon/2]))\munion(K_t\minter B_r(p)^c),
$$
\noindent then $K_0\subseteq K_t'$. However, $K_t'$ is a compact, convex subset of $\hat{K}$. Moreover, being locally the intersection of two convex sets of type $F_t$, by Lemma \procref{LemmaIntersectionIsOfTypeF}, $K_t$ is also of type $F_t$. In particular, it is an element of $\Cal{F}$ which is a strict subset of $K_t$, which yields a contradiction. The result follows.\qed
\medskip
\noindent We thus obtain Theorem \procref{TheoremExistence}:
\medskip
{\noindent\bf Proof of Theorem \procref{TheoremExistence}:\ }Lemma \procref{LemmaExistence} yields convex sets whose boundaries are of constant Gaussian curvature in the viscosity sense (c.f. \cite{HarveyLawsonI} and \cite{HarveyLawsonII}). By Theorem $5.1$ of \cite{GuanSpruckI}, the boundaries of these sets are smooth, and this completes the proof.\qed
\goodbreak
\newhead{Bibliography}
{\leftskip = 5ex \parindent = -5ex
\leavevmode\hbox to 4ex{\hfil \cite{HarveyLawsonI}}\hskip 1ex{Harvey F. R., Lawson H. B. Jr., Dirichlet Duality and the Nonlinear Dirichlet Problem, {\sl Comm. Pure Appl. Math.} {\bf 62} (2009), no. 3, 396--443}
\medskip
\leavevmode\hbox to 4ex{\hfil \cite{HarveyLawsonII}}\hskip 1ex{Harvey F. R., Lawson H. B. Jr., Dirichlet Duality and the Nonlinear Dirichlet Problem on Riemannian Manifolds, arXiv:0907.1981}
\medskip
\leavevmode\hbox to 4ex{\hfil \cite{CaffNirSprV}}\hskip 1ex{Caffarelli L., Nirenberg L., Spruck J., Nonlinear second-order elliptic equations. V. The Dirichlet problem for Weingarten hypersurfaces, {\sl Comm. Pure Appl. Math.} {\bf 41} (1988), no. 1, 47--70}
\medskip
\leavevmode\hbox to 4ex{\hfil \cite{CrandIshiiLions}}\hskip 1ex{Crandall M. G., Ishii H., Lions P.-L., User's guide to viscosity solutions of second order partial differential equations, {\sl Bull. Amer. Math. Soc.} {\bf 27} (1992), no. 1, 1--67}
\medskip
\leavevmode\hbox to 4ex{\hfil \cite{GuanSpruckI}}\hskip 1ex{Guan B., Spruck J., The existence of hypersurfaces of constant Gauss curvature with prescribed boundary, {\sl J. Differential Geom.} {\bf 62} (2002), no. 2, 259--287}
\medskip
\leavevmode\hbox to 4ex{\hfil \cite{GuanSpruckII}}\hskip 1ex{Guan B., Spruck J., Locally Convex Hypersurfaces of Constant Curvature with\break Boundary, {\sl Comm. Pure Appl. Math.} {\bf 57} (2004), 1311--1331}
\medskip
\leavevmode\hbox to 4ex{\hfil \cite{SmiCGC}}\hskip 1ex{Smith G., Constant Gaussian Curvature Hypersurfaces in Hadamard Manifolds,\break arXiv:0912.0248}
\medskip
\leavevmode\hbox to 4ex{\hfil \cite{SmiSLC}}\hskip 1ex{Smith G., Special Lagrangian curvature, arXiv:math/0506230}
\medskip
\leavevmode\hbox to 4ex{\hfil \cite{SmiNLD}}\hskip 1ex{Smith G., The Non-Linear Dirichlet Problem in Hadamard Manifolds,\break arXiv:0908.3590}
\par}
\enddocument

%% file: references.tex
\global\def\_@citation@HarveyLawsonI{1}
\global\def\_@citation@HarveyLawsonII{2}
\global\def\_@citation@CaffNirSprV{3}
\global\def\_@citation@CrandIshiiLions{4}
\global\def\_@citation@GuanSpruckI{5}
\global\def\_@citation@GuanSpruckII{6}
\global\def\_@citation@SmiCGC{7}
\global\def\_@citation@SmiSLC{8}
\global\def\_@citation@SmiNLD{9}
\global\def\_@proc@TheoremExistence{1.1}
\global\def\_@head@HeadPSets{2}
\global\def\_@proc@LemmaCharacterisation{3.1}
\global\def\_@proc@LemmaIntersectionIsOfTypeF{3.3}
\global\def\_@proc@LemmaExistence{4.1}